\documentclass{amsart}

\makeindex

\usepackage[all]{xy}
\usepackage{amssymb}

\usepackage{amsmath}

\usepackage{amsthm}

\usepackage{amscd}

\usepackage{amsfonts}

\usepackage{amssymb}

\usepackage{mathrsfs}

\newtheorem{Lem}{Lemma}[section]

    \newtheorem{Prop}[Lem]{Proposition}

    \newtheorem{Thm}[Lem]{Theorem}

    \newtheorem{Cor}[Lem]{Corollary}

\theoremstyle{definition}

    \newtheorem{Rem}[Lem]{Remark}

\newcommand{\Ps}{\mathbb{P}^}
\newcommand{\Z}{\mathbb{Z}}

\newcommand{\ra}{\rightarrow}
\newcommand{\lra}{\longrightarrow}
\newcommand{\R}{\mathbb R}

\newcommand{\C}{\mathbb C}

\begin{document}
\title{A one dimensional family of $K3$ surfaces\\ with a $\Z_4$ action}
\author{Michela Artebani}
\address{Dipartimento di Matematica, Universit\`a degli Studi di Milano, via C. Saldini 50, Milano, Italia}
\email{michela.artebani@unimi.it} \subjclass{} \keywords{genus three
curves, K3 surfaces}
\thanks{This work is partially supported by: 
PRIN 2005: {\em Spazi di moduli e
teoria di Lie};  Indam (GNSAGA) and the NSERC Discovery Grant of Noriko Yui
at Queen's University, Canada.}
\date{\today}
\begin{abstract}
The minimal resolution of the degree four cyclic cover of the plane
branched along a GIT stable quartic is a K3 surface with a non symplectic action of $\Z_4$. In this paper
we study the geometry of the one dimensional family of K3 surfaces
associated to the locus of plane quartics with five nodes.
\end{abstract}
\maketitle 
\pagestyle{myheadings} 
\markboth{Michela Artebani}{A one
dimensional family of $K3$ surfaces with a $\Z_4$ action}
\setcounter{tocdepth}{1}
\section*{Introduction}
Let $V\subset\mid\!\mathcal O_{\Ps 2}(4)\!\mid$ be the space of
plane quartics with five nodes and $\mathcal V$ be the one dimensional family given by the quotient of $V$ for the action of $PGL(3,\C)$.
The minimal resolution $X_C$ of the degree four cyclic cover of the
plane branched along a quartic $C\in V$ is a $K3$ surface equipped with a non-symplectic automorphism group $G_C\cong\Z_4$. Since the isomorphism class of $X_C$ only depends on the projective
equivalence class of $C$, this construction gives a
one dimensional family $\mathcal X$ of $K3$ surfaces. Moreover, as proved in \cite{m}, it defines an injective period map
$$P:\mathcal V\lra \mathcal M$$
where $\mathcal M$ is a moduli space for couples $(X_C,G_C)$
(as defined in \cite{DK}).

This paper describes the geometry of the family $\mathcal X$
by studying the structure of the moduli space $\mathcal M$, the behavior of the period map on the closure of $\mathcal V$ and the occurence of singular $K3$ surfaces.

In the first section we introduce the $K3$ surface associated to a GIT stable plane quartic according to the construction given by S. Kond\=o in \cite{K}.

The period domain of these polarized $K3$ surfaces is isomorphic to the complex one dimensional ball, the second section shows that their moduli space $\mathcal M$ is the Fricke modular curve of level two.

Any $K3$ surface $X_C$ carries an elliptic fibration induced by the pencil of lines through one node of $C\in V$. In section $3$ we prove that the fibration is isotrivial and the generic fiber is isomorphic to the elliptic curve $E=\C/\Z[i]$. In fact, after a base change and a normalization, the fibration is the product $B_C\times E$ where $B_C$ is a genus two curve with splitting Jacobian $J(B_C)=E_C\times
E_C$. 

In section $4$ we describe the behavior of the period map on the
closure $\overline{\mathcal V}$ of $\mathcal V$.
We prove that the period map can be extended to $\overline{\mathcal V}$,
giving an isomorphism with the projective line.

The last section shows that there is a correspondence between  $X_C$ and the Kummer surface $Km(E\times E_C)$. In particular, the occurence of singular $K3$ surfaces in the family is due to the existence of isogenies between the elliptic curves $E$ and $E_C$. This is also
connected to the existence of certain rational ``splitting curves'' for $C$ (see \cite{H}).
We finally give a partial characterization of transcendental lattices of singular $K3$ surfaces in the family. In particular we prove that the Fermat quartic, the Klein quartic and Vinberg's K3 surface (see \cite{V}) belong to the family.\\

\emph{Acknowledgements.}  
I'm grateful to B. van Geemen for introducing me to the subject and for many interesting discussions. I also would like to thank I. Dolgachev, A. Laface and N. Yui for several helpful comments.

\section{Plane quartics and $K3$ surfaces}\label{def}
Let $C$ be a \emph{GIT stable} plane quartic i.e. having at most ordinary nodes and cusps (see \cite{Gi}). The degree four cyclic cover of the plane branched along $C$ has at most rational double points, hence its minimal resolution $X_C$ is a $K3$ surface with an order four automorphism group $G_C$ (see \cite{K} and \cite{m}). In this
paper we consider the locus $V$ of stable plane quartics with five
nodes. By taking the quotient of $V$ for the natural action of
$PGL(3,\C)$ we get a one dimensional family $\mathcal V$. The
isomorphism class of the cover only depends on the class of the
quartic in $\mathcal V$, hence this construction defines a period map$$P: \mathcal V\lra \mathcal M,\ \ \left[ C \right] \longmapsto \left[(X_C,G_C)\right]$$
where $\mathcal M$ is a one dimensional moduli space parametrizing couples $(X_C,G_C)$ (the precise definition is given in \ref{mod}).

We now choose a parametrization for $\mathcal V$.
Consider the plane quartic $C_{\alpha}$, $\alpha\in\Ps 1$ which is the union of the following conic and two lines:
$$
Q:\ y^2-xz=0,\  \ L:\  y=0,\ \ M_{\alpha}:\ \alpha x+2y+z=0.$$
This gives a one parameter non-constant family of plane quartics with five nodes, hence the general point in $\mathcal V$ is represented by a curve in this family.

We denote with $\pi_{\alpha}$ the four cyclic cover of the plane branched along $C_{\alpha}$
$$\pi_{\alpha}:Y_{\alpha}\lra \Ps 2$$ and with  $\nu_{\alpha}$ its minimal resolution
$$\nu_{\alpha}:X_{\alpha}\lra Y_{\alpha}.$$
By the previous remark, the general $X_{\alpha}$ is a $K3$ surface.
Let $G_{\alpha}\cong \Z_4$ be the order four automorphism group of covering transformations of $\pi_{\alpha}$.




\section{The period domain}
In this section we describe the moduli space parametrizing couples $(X_{\alpha},G_{\alpha})$ where $X_{\alpha}$ is a $K3$ surface associated to a plane quartic with five nodes and $G_{\alpha}$ is the corresponding order four covering transformation group.

Let $\sigma_{\alpha}$ be a generator of $G_{\alpha}$ and $\rho_{\alpha}$ be the induced isometry on the cohomology lattice $H^2(X_{\alpha},\Z)$.
In \cite{m} it is proved that $\sigma_{\alpha}$ acts as a primitive 4-th root of unity on $H^{2,0}(X_{\alpha})$. In fact we can assume
$$\rho_{\alpha}(\omega_{\alpha})=i\omega_{\alpha}$$
where $\omega_{\alpha}\not=0$ is a holomorphic two-form on $X_{\alpha}$.
In particular, the invariant lattice of the involution $\tau_{\alpha}=\sigma_{\alpha}^2$ is contained in the Picard lattice of $X_{\alpha}$. In fact we will show that the invariant lattice is the Picard lattice of the generic $K3$ surface $X_{\alpha}$.
\subsection{The generic point}
Let $T_{\alpha}$ and $N_{\alpha}$ be the
transcendental lattice and the Picard lattice of the generic
$K3$ surface $X_{\alpha}$ respectively.
\begin{Lem}\label{iso}
The isomorphism classes of $T_{\alpha}$ and $N_{\alpha}$ are given by
$$T_{\alpha}=A_1^{\oplus 2}\oplus A_1(-1)^{\oplus 2},\
N_{\alpha}=U\oplus E_7^{\oplus 2}\oplus A_1^{\oplus 2}.$$
Moreover, in the natural basis of $T_{\alpha}$, the action of the
isometry $\rho_{\alpha}$ is given by the matrix:
$$J=A\oplus A$$
where:
$$A=\left(\begin{array}{cc}
0&-1\\
1 &0\\
\end{array}\right).$$
\end{Lem}
\proof
Let $\iota_{\alpha}=\rho_{\alpha}^2$ and $L_{\pm}(\alpha)\subset
H^2(X_\alpha,\Z)$ be its eigenspaces. It can be easily seen that the
fixed locus of the involution $\tau_{\alpha}$ is the disjoint union of eight smooth rational curves. By Theorem 4.2.2, \cite{N2} this implies:
$$r(L_{+}(\alpha))=18,\ \ell(L_{+}(\alpha))=4,$$
where $r(\cdot)$ denotes the rank of a lattice and $\ell(\cdot)$ the minimal number of generators of its discriminant group.
Since the family $\{C_{\alpha}\}_{\alpha\in \Ps 1}$ is one
dimensional, the Picard lattice equals $L_+(\alpha)$ for
general $\alpha$, in particular $r(T_{\alpha})=4$. By Theorem 3.1, \cite{N3} there is an isomorphism of $\Z[i]$-modules:
$$T_{\alpha}\cong \Z[i]\oplus \Z[i].$$
Notice that an even symmetric lattice $\Lambda$
which is a free $\Z[i]$-module of rank one with $i\in O(\Lambda)$,
$i^2=-id$ is of the form:
$$\Lambda\cong A_1(n)\oplus A_1(n),\ \ n\in\Z,$$
where the action of the isometry $i$ is given by the matrix $A$.
Hence, in a suitable integral basis, the transcendental lattice
$T_{\alpha}$ has intersection matrix of the form:
$$B=\left(\begin{array}{cc}
A_1(n)^{\oplus 2}&C\\
C& A_1(m)^{\oplus 2}
\end{array}\right)
,$$ where:
$$C=\left(\begin{array}{cc}
b & c\\
-c &b\\
\end{array}\right)
.$$
Since $T_{\alpha}$ is a 2-elementary lattice with
$\ell(T_{\alpha})=4$ we have $det(B)=2^4$. Moreover, its signature
is $(2,2)$. This implies $b=c=0$ and $nm=-1$. Hence $T_{\alpha}\cong
A_1^{\oplus 2}\oplus A_1(-1)^{\oplus 2}$. In particular, the Picard lattice is an
even hyperbolic 2-elementary lattice with the invariants
$(s_+,s_-,\ell,\delta)=(1,17,4,1)$. By Theorem 3.6.2, \cite{N1} it is
isomorphic to the lattice $U\oplus E_7^{\oplus 2}\oplus A_1^{\oplus
2}$.\qed
\subsection{The moduli space}\label{mod}
Let $L$ be the abstract $K3$ lattice and $\rho=\rho_{\alpha}$, $\tau=\tau_{\alpha}$. We denote with $N$ and $T$ the positive and negative eigenlattices of $\tau$ on $L$ respectively.
By the remarks in the previous section it follows that the period domain for $K3$ surfaces in the family is given by:
$$D=\{z\in \mathbb P(T\otimes C): \rho(z)=iz, (z,\bar z)>0\}.$$
Since $T$ has rank $4$ it can be easily seen that $D$ is a one dimensional complex ball.
By taking the quotient for the arithmetic group:
$$\Gamma=\{\gamma\in O(T): \gamma\circ\rho=\rho\circ\gamma\}$$
we get the moduli space
$$\mathcal M=D/\Gamma.$$
Let $T_{-2}=\{\delta\in T: \delta^2=-2\}$, $H_{\delta}=\delta^{\perp}\cap D$ and
$$\Delta=\bigcup_{\delta\in T_{-2}}H_{\delta}.$$
\begin{Prop}
The quotient $(D\backslash \Delta)/\Gamma$ parametrizes isomorphism classes of couples $(X_{\alpha},G_{\alpha})$.
Moreover, the period map $P:\mathcal V\lra (D\backslash \Delta)/\Gamma$ is an isomorphism.
\end{Prop}
\proof See \cite{m}.\qed\\

Lemma \ref{iso} allows us to describe in detail the structure of the moduli space $\mathcal M$.
Consider the following subgroups of $SL(2,\mathbb Z)$:
$$G_0=SU(1,1)\cap M(2,\mathbb Z[i]),$$
$$H_0=\{\left(\begin{array}{cc}
a & b\\
c &d
\end{array}\right)\in SL(2,\mathbb Z): a+d\equiv b+c \equiv 0\ (mod\ 2)\}.$$
\begin{Prop}
We have the isomorphisms:
$$\mathcal M\cong B/G\cong S/H$$
where $B=\{z\in \mathbb C: \mid z\mid<1\}$ is the complex 1-ball, $S=\{z\in \mathbb C:\ Im(z)>0\}$ is
the Siegel upper half space and
$$G=G_0\cup LH_0,$$
$$H=H_0\cup MH_0$$
where
$$L=\left(\begin{array}{cc}
e^{-i\pi/4} & 0\\
0 & e^{i\pi/4}
\end{array}\right),$$
$$M=\frac{1}{\sqrt 2}\left(\begin{array}{cc}
1& -1\\
1 &1
\end{array}\right).$$
\end{Prop}
\proof The period domain $D$ is given by points
$z=(z_1,\dots,z_4)\in \mathbb P(T\otimes \mathbb C)$ such
that:
$$1)\ ^{t}z T\bar z>0$$
$$2)\ Jz=iz.$$
Hence $z$ is of the form:
$$z=(iz_2,z_2,iz_4,z_4),\ \ \mid z_2\mid ^2-\mid z_4\mid^2>0.$$
Thus we get the isomorphism:
$$\Psi_1:D\rightarrow B=\{w\in\mathbb C:\mid w\mid <1\},\ \ z\mapsto z_4/z_2.$$
We are interested in the following subgroup of the isometries of
$T$:
$$\Gamma=\{M\in O(T): MJ=JM\}.$$
Under the identification:
$$\mathbb Z[J]\cong \mathbb Z[i]$$
we have the isomorphism $T\cong \mathbb Z[i]^{2}$ as
$\mathbb Z[i]$-modules. It can be easily seen that in the natural
basis for $\mathbb Z[i]^2$ the intersection form on $T$ is
given by:
$$Q(z,w)=2(z\bar z-w\bar w).$$
Then we get:
$$\Gamma=U(Q)\cap M_2(\mathbb Z[i])\cong U(1,1)\cap M_2(\mathbb Z[i]).$$
Let $M\in \Gamma$:
$$M=\left(\begin{array}{cc}
a & b\\
c &d
\end{array}\right)$$
where $a,b,c,d\in \mathbb Z[i]$. The action of $M$ on $D$
induces an action of $M$ on $B$ which is given by the
M\"{o}bius transformation:
$$z\mapsto \psi_M(z)=\frac{c+dz}{a+bz}.$$
Since two matrices in $M_2(\mathbb C)$ give the same M\"{o}bius
transformation if and only if they are the same up to multiplication
for a nonnegative scalar,  the group of M\"obius transformations of
$\mathbb C$ is isomorphic to the quotient:
$$\mathcal T\cong SL(2,\mathbb C)/{\pm I}.$$
Consider the homomorphism:
$$T:GL(2,\mathbb C) \rightarrow \mathcal T \ \ \ M\mapsto \frac{1}{\sqrt{det(M)}}M.$$
Notice that the kernel of $T$ is isomorphic to $\mathbb C^*$. Let
$T_{\mid\Gamma}$ be the restriction of $T$ to $\Gamma$, then
$ker(T_{\mid \Gamma})\cong F_4$ where $F_4$ is isomorphic to the
group of 4-th roots of unity. Notice that
$G=Im(T_{\mid\Gamma})\subset SU(1,1)/{\pm I}$ is given by:
$$G=\{M\in SU(1,1)/{\pm I}\mid \exists\epsilon \in \mathbb C^*: \epsilon M \in M(2,\mathbb Z[i])\}.$$
Let:
$${\Gamma}_0=SU(1,1)\cap M(2,\mathbb Z[i])\subset G$$
and $G_0$ its image in $\mathcal T$. Let:
$$L'=\left(\begin{array}{cc}
1 & 0\\
0 &i
\end{array}\right)\in \Gamma$$
and
$$L=T(L')=[ \left(\begin{array}{cc}
e^{-i\pi/4} & 0\\
0 & e^{i\pi/4}
\end{array}\right)].$$
Notice that ${L}^{-1}{G_0}L=G_0$ and ${L}^2\in G_0$. In fact, if
$M\in G$ then $LM\in G_0$. Hence:
$$G=G_0\cup L G_0.$$

A biholomorphic map between $B$ and $S=\{z\in
\mathbb C: Im(z)>0\}$ is given by the M\"obius transformation
$\Psi_2=\psi_K$:
$$\psi_K:B\rightarrow S\ \ \ z\mapsto \frac{i+z}{1+iz}.$$
associated to the matrix:
$$K=\left(\begin{array}{cc}
 1& i\\
i & 1
\end{array}
\right) .$$ The map $\psi_K$ induces the isomorphism between the
groups of automorphisms:
$$\Upsilon:Aut(B)\rightarrow Aut(S)\ \  \ \phi\mapsto \phi'=\psi_K \phi {\psi_K}^{-1}.$$
Let $\psi_M$ be the M\"obius transformation corresponding to a
matrix $M\in SU(1,1)$:
$$M=\left(\begin{array}{cc}
a & b\\
\bar b &\bar a
\end{array}\right)$$
where $a,b\in \mathbb C$. Then the map $\phi'$ is the M\"{o}bius
transformation associated to the matrix $KMK^{-1}\in SL_2(\mathbb
\mathbb R)$:
$$KMK^{-1}=\left( \begin{array}{cc}
Re(a)+Im(b) & Re(b)+Im(a)\\
Re(b)-Im(a) & Re(a)-Im(b)
\end{array} \right).$$
Conversely, let $N\in SL_2(\mathbb R)$:
$$N=\left( \begin{array}{cc} \alpha &\beta\\
\gamma & \delta
\end{array}\right).
$$
Then $N=KMK^{-1}$ where $M\in SU(1,1)$ is given by:
$$a=\frac{1}{2}[\alpha+\delta+i(\beta-\gamma)],\ b=\frac{1}{2}[\beta+\gamma+i(\alpha-\delta)].$$
This gives an isomorphism between the groups of M\"obius
transformations associated to $SU(1,1)$ and that associated to
$SL(2,\mathbb R)$.
The image of $G_0$ is the following subgroup of $H$:
$$H_0=\{\left(\begin{array}{cc}
\alpha & \beta\\
\gamma &\delta
\end{array}\right)
\in SL(2,\mathbb Z): \alpha+\delta\cong \beta+\gamma\cong 0\ (mod\
2)\}.$$ The image of $L$ in $SL(2,\mathbb R)$ is given by:
$$\Upsilon(L)=\frac{1}{\sqrt 2}\left(\begin{array}{cc}
1 & -1\\
1 & 1
\end{array}\right).$$
Then we have the following description:
$$H=H_0\cup \Upsilon(L)H_0.$$
\qed\\
Consider the level $2$ congruence subgroup:
$$H[2]=\{\left(\begin{array}{cc}
a & b\\
c &d
\end{array}\right)\in PSL(2,\mathbb Z)\mid c\cong 0\ (mod\ 2)\}.$$
The order two element
$$F=\left(\begin{array}{cc}
0&-1/\sqrt{2}\\
\sqrt{2}& 0
\end{array}\right)\in PSL(2,\R)$$
lies in the normalizer of $H[2]$ in $PSL(2,\R)$ and it is called
\emph{Fricke involution}. The group:
$$H[2]^+=H[2]\cup FH[2]\subset PSL(2,\R)$$
is called \emph{Fricke modular group of level 2} and the
quotient:
$$C(2)^+=\mathcal S/H[2]^+$$
is the \emph{Fricke modular curve of level 2}.
\begin{Cor}
We have the isomorphisms:
$$\mathcal M\cong C(2)^+\cong \mathbb A^1.$$
\end{Cor}
\proof The group $H_0$ is conjugated to $H[2]$ in $SL(2,\mathbb Z)$:
$$TH_0T^{-1}=H_0\ \ ,T=\left(\begin{array}{cc}
1 & 0\\
1 &1
\end{array}\right).$$
Besides, it can be easily proved that:
$$T(\Upsilon(L)H_0)T^{-1}=FH_0.$$
Hence the group $H$ is isomorphic to the Fricke modular group of
level two and $\mathcal M$ is isomorphic to $C(2)^+$. The last isomorphism follows by
Proposition 7.3 and Corollary 7.4, \cite{D}.
\qed
\begin{Rem}
In \cite{D} it is proved that the Fricke modular curve of level $2$ is also the moduli space for the mirror family of degree $4$ polarized $K3$ surfaces. It would be interesting to understand if there is any geometric correspondence between the two families.
\end{Rem}
 \section{An elliptic pencil}\label{pencil}
In this section we show that $X_{\alpha}$ carries a natural elliptic fibration, induced by the pencil of lines through one of the nodes of $C_{\alpha}$.
\subsection{Definition}
Note that the conic $Q$ intersects the line $L$ in
$p_1=(0:0:1)$ and $p_2=(1:0:0)$.
\begin{Prop}\label{fi}
The pencil of lines through the point $p_1$ induces an isotrivial
elliptic fibration $\mathcal E_{\alpha}: X_{\alpha}\lra \Ps 1$. After a base change $B_{\alpha}\lra \Ps 1$ and a
normalization, the fibration is the trivial fibration:
$$E\times B_{\alpha}\lra B_{\alpha},$$
where $E=\C/\Z[i]$ is the elliptic curve with $j=1728$ and
$B_{\alpha}$ is a genus two curve.
\end{Prop}
\proof
The pencil of lines through the point $p_1\in Q\cap L$ is given by
the equation:
$$
y=\lambda x,
$$
where $\lambda\in \Ps1$.
We substitute $y=\lambda x$ in the equation of $C_{\alpha}$
and we restrict to the affine subset where $x=1$:
$$
\lambda(z-\lambda^2)(z+(2\lambda+\alpha))=0.
$$
In general we have:
$$
(z-a)(z-b)=z^2-(a+b)z+ab= (z-\mbox{$\frac{1}{2}$}(a+b))^2-
(\mbox{$\frac{1}{2}$}(a-b))^2 .$$ Introducing a new variable $z_1$
by:
$$
z=\mbox{$\frac{1}{2}$}(a-b)z_1+\mbox{$\frac{1}{2}$}(a+b) ,$$ we get:
$$
(z-a)(z-b)=\mbox{$\frac{1}{4}$}(a-b)^2\,(z_1^2-1).
$$
In our case $a=\lambda^2,b=-(2\lambda+\alpha)$, so:
$$
\lambda(z-\lambda^2)(z+(2\lambda+\alpha))=
\mbox{$\frac{1}{4}$}\lambda(\lambda^2+2\lambda+\alpha)^2(z_1^2-1).
$$
Thus we are considering the the fibration on $Y_{\alpha}$:
$$
Y_{\alpha}\lra \Ps 1_\lambda,\qquad
w^4=\mbox{$\frac{1}{4}$}\lambda(\lambda^2+2\lambda+\alpha)^2(z_1^2-1).
$$
This induces an elliptic fibration on $X_{\alpha}$ with fibers
isomorphic to the elliptic curve:
$$E:\ w^4=z_1^2-1.$$
Notice that $E\cong \C/\Z[i]$ since $E$ has an automorphism of order
$4$ which fixes a point:
$$(z_1,w)\longmapsto(z_1,iw),$$
(the point $(z_1,w)=(1,0)$ is fixed). 

To get the trivial fibration, we first make a base change
$$
\Ps1_\rho\lra \Ps1_\lambda,\ \  \rho\longmapsto \lambda=\rho^2 ,$$
which gives the equation:
$$
w^4=\left(\mbox{$\frac{1}{2}$}\rho(\rho^4+2\rho^2+\alpha)\right)^2
(z_1^2-1).
$$
Next we consider the genus two curve:
$$
B_\alpha:\ \tau^2=\rho(\rho^4+2\rho^2+\alpha)).
$$
We make the base change:
$$
B_\alpha\longrightarrow \mathbb P^1_\rho,\ \  (\rho,\tau)\longmapsto
\rho
$$
and we define $w=\tau w_1/\sqrt{2}$, so we get
$$
w_1^4=z_1^2-1.
$$
Hence the normalization of the pull-back of the family is the
product $B_\alpha\times E$. \qed

\subsection{The Weierstrass model}\label{wfbis}
We now determine the Weierstrass model for the isotrivial elliptic fibration
defined in Proposition \ref{fi}:
\begin{Lem}
The Weierstrass form for the elliptic fibration $\mathcal E_{\alpha}$ is
given by:
$$
Y_{\alpha}\lra \Ps1_\lambda,\ \
v^2=u^3-\lambda^3(\lambda^2+2\lambda+\alpha)^2 u.
$$
\end{Lem}
\proof
Recall that $\mathcal E_{\alpha}$ is given by
$$\mathcal E_{\alpha}:Y_{\alpha}\lra\Ps1_\lambda, \ \
w^4=\mbox{$\frac{1}{4}$}\lambda(\lambda^2+2\lambda+\alpha)^2(z_1^2-1).$$
Let
$$
\beta=\left(\mbox{$\frac{1}{4}$}\lambda(\lambda^2+2\lambda+\alpha)\right)^{-2},
$$
then the fibration can be rewritten as:
$$
z_1^2=\beta w^4+1.
$$
Introducing the new coordinates $w=1/s$, $z_1=t/s^2$  we get:
$$
t^2=s^4+\beta\Longrightarrow(t-s^2)(t+s^2)=\beta.
$$
Next we put (see \cite{Ca}):
$$
x=t+s^2\Longrightarrow t-s^2=\beta/x,\  2s^2=x-\beta/x.
$$
Multiply the last equation by $x^2$ and put $y=sx$:
$$
2v^2x^2=x^3-\beta x\quad\Longrightarrow\quad 2y^2=x^3-\beta x.
$$
We finally put $x=u/2,\;y=v/4$ and multiply the equation by $8$:
$$
v^2=u^3-4\beta u,
$$
where we have:
$$
u=2 \frac{z_1+1}{w^2},\  v=4 \frac{z_1+1}{w^3}.
$$
Hence the family is:
$$
v^2=u^3-\frac{16}{\lambda(\lambda^2+2\lambda+\alpha)^2}u .$$
The transformation to the form
$w^4=\mbox{$\frac{1}{4}$}\lambda(\lambda^2+2\lambda+\alpha)^2(z_1^2-1)$
is given by:
$$
u=\left(\frac{\lambda(\lambda^2+2\lambda+\alpha)}{w}\right)^2\frac{z_1+1}{2},\
v=\left(\frac{\lambda(\lambda^2+2\lambda+\alpha)}{w}\right)^3\frac{z_1+1}{2}.
$$
The original variables $x,y,z,w$ can be obtained from:
$$
z_1=2z-\frac{\lambda^2-2\lambda-\alpha}{\lambda^2+2\lambda+\alpha}.
$$
\qed
\subsection{Singular fibers}
The singular fibers in a Weierstrass fibration with equation:
$$v^2=u^3-f(\lambda)u,\ \ f\in\C[\lambda],$$
correspond to the values $\lambda$ where $f(\lambda)=0$ and to
$\lambda=\infty$ if $\deg(f)$ is not divisible by $4$.  Let
$f=(\lambda-a)^kg$ with $g(a)\neq 0$, then we may always assume that
$0\leq k\leq 3$ and we have bad reduction in $a$ only if $k\neq 0$.
Thus we get three types of bad fibers for $k=1,2,3$, they are:
\begin{itemize}
\item{} $k=1$, type $III$ (two tangent rational curves), $\chi=3$,
\item{} $k=2$, type $I_0^*$ (a double component and 4 reduced comp.), $\chi=6$,
\item{} $k=3$, type $III^*$ (eight components), $\chi=9$.
\end{itemize}
It is now easy to find the bad fibers in our case:
$$v^2=u^3-\lambda^3(\lambda^2+2\lambda+\alpha)^2u.$$
\begin{Lem}\label{sing}
The elliptic fibration $\mathcal E_{\alpha}$ has the following configuration of singular fibers:
\begin{itemize}
\item{} type $III^*$ over $\lambda=0$,
\item{} type $I_0^*$ over the solutions $\lambda_1,\lambda_2$ of $\lambda^2+2\lambda+\alpha$,
\item{} type $III$ over $\lambda=\infty$.
\end{itemize}
\end{Lem}
Notice that, for $\lambda=0$ we get the line $L$, for
$\lambda=\infty$ we get the line $x=0$, which is tangent to the
conic $Q$ in $p_1$. Finally, the values $\lambda_i$, $i=1,2$
correspond to lines through the intersection points of $Q$ and
$M_{\alpha}$.
\begin{Rem}\ \\
i) The elliptic fibration $\mathcal E_{\alpha}$ has two sections,  given by the line $M_{\alpha}$ and the conic $Q$, which cut on each fiber the two fixed points of the order $4$ automorphism (defined by $(z_1,w)=(\pm 1,0)$).\\
ii) It follows easily from Lemma \ref{iso} that the generic $K3$ surface $X_{\alpha}$ has also an elliptic fibration with two fibers of type $III^*$ and two of type $III$. In fact, this is the elliptic fibration induced by the pencil of lines through the intersection point of $L$ and $M_{\alpha}$.  \end{Rem}
\section{Compactification}
If the quartic $C_{\alpha}$ is not stable, then $X_{\alpha}$ is not a $K3$ surface. However, we show that in some cases proper modifications of the family still give $K3$ surfaces in the limit.
In other words, we study the behavior of the period map $P$ on the closure $\overline{\mathcal V}$ of $\mathcal V$.

Note that $M_\alpha$ (see section \ref{def}) is the pencil of lines through the point $(0,1,-2)$ and the curve $C_{\alpha}$ is not stable iff $M_\alpha$ is tangent to $Q$ or if it contains a point of $Q\cap L$.
The line $M_{\infty}$ is tangent to $Q$ in the point $p_1\in Q\cap
L$. Furthermore, the line $M_1$ is tangent to $Q$ in the point
$(1,-1,1)$, hence $C_1$ has a tacnode in this point. Finally, the
line $M_0$ passes through the point $(1,0,0)\in Q\cap L$, hence
$C_0$ has a triple point. Thus there are three unstable quartics in the family:
$$
\alpha=0,1,\infty.
$$


\begin{Lem}\label{inf}
There exists a modification $X'_{\alpha}$ of the family $X_{\alpha}$ such that the fiber $X'_{\infty}$ is a $K3$ surface with \\
i) an elliptic fibration with the same configuration of singular fibers of Lemma \ref{sing}; \\
ii) an automorphism of order eight with transcendental value $m=8$.\\
iii) Picard number 18.
\end{Lem}
\proof We consider the elliptic fibration in Weierstrass form from
section \ref{wfbis}
$$
\mathcal E_{\alpha}:Y_{\alpha}\longrightarrow \mathbb P^1_\lambda,\qquad
v^2=u^3-\lambda^3(\lambda^2+2\lambda+\alpha)^2 u.
$$
We put
$$
\alpha:=\beta^{-8},\qquad u:=\beta^{-14}u,\quad v=\beta^{-21}v,\quad
\lambda=\beta^{-4}\lambda.
$$
Then, after multiplying the equation by $\beta^{-42}$, we get:
$$
Y'_\beta:\quad v^2=u^3-\lambda^3 (\lambda^2+2\beta^4\lambda+1)^2u.
$$
This modified family has a good reduction for $\beta\rightarrow 0$.
The fibration $Y'_{\infty}\lra \mathbb P^1_\lambda$ has 4 bad
fibers with the same configuration of the general case.
Moreover, the surface $Y'_{\infty}$ has an extra automorphism $\varphi$ given
by:
$$
u:=\zeta^2 u, \qquad v:=\zeta^3 v,\qquad \lambda:= -\lambda\qquad
(\zeta^4=-1).
$$
Note that the holomorphic two form on $X'_{\infty}$ is locally given
by $\omega=(d\lambda\wedge du)/v$ and
$\varphi^*\omega=(-1\zeta^2/\zeta^3)\omega=-\zeta^{-1}\omega$. Hence
the transcendental value is equal to eight and the transcendental
lattice of $X'_{\infty}$ allows the action of the ring
$\Z[\zeta]$ (see Theorem 3.1, \cite{N3}), in particular its rank is a multiple of 4. \qed
\begin{Lem}\label{0}
There exists a modification $X''_{\alpha}$ of the family $X_{\alpha}$
such that the fiber $X''_{0}$ is the ``Vinberg's K3 surface'' and carries an elliptic
fibration with two fibers of type $III^*$ and one of type $I_0^*$.
\end{Lem}
\proof We consider again the elliptic fibration in Weierstrass form
from section \ref{wfbis}
$$
\mathcal E_{\alpha}\longrightarrow \mathbb P^1_\lambda,\qquad
v^2=u^3-\lambda^3(\lambda^2+2\lambda+\alpha)^2 u.
$$
When $\alpha\rightarrow 0$ we get $\lambda^5(\lambda+2)$, and
changing coordinates allows to reduce to the case
$\lambda(\lambda+2)$, which gives no longer a $K3$ surface.
We consider the fibration near $\lambda=\infty$, so we put:
$$
\mu=\lambda^{-1},\quad u:=u/\mu^4,\quad v:=v/\mu^6
$$
and multiply throughout by $\mu^{12}$:
$$
v^2=u^3-\mu(1+2\mu+\alpha\mu^2)^2 u.
$$
We make a base change and a coordinate change:
$$
\alpha=\beta^4,\qquad \mu:=\mu/\beta^4,\quad u:=u/\beta^6,\quad
y:=y/\beta^9,
$$
and multiply throughout by $\beta^{18}$:
$$
v^2=u^3-\mu(\beta^4+2\mu+\mu^2)^2 u.
$$
It is now obvious that for $\beta\rightarrow 0$ we get an elliptic
fibration on a $K3$ surface $X''_0$ associated to:
$$v^2=u^3-\mu^3(2+\mu)^2 u.$$
Notice that there are $2$ fibers of type $III^*$ over $\mu=0,\infty$ and one of type $I_0^*$ over $\mu=-1/2$.
It follows from the Shioda-Tate formula (Corollary 1.5, \cite{Shi}) that the rank of the Picard lattice of $X''_0$ is at least $20$. The table of Shimada-Zhang (see \cite{SZ}, case 279) shows that a fibration with this fiber type is unique and has
transcendental lattice isomorphic to $A_1(-1)^{\oplus 2}$. Hence,
the surface $X''_0$ is the ``Vinberg's K3 surface'' (see \cite{V}).

\begin{Prop}
The period map $P$ can be extended to an isomorphism
$$P:\overline{\mathcal V}\lra \overline{\mathcal M}$$
where $\overline{\mathcal M}=\mathcal M\cup\{X\}\cong \Ps 1$ is the Baily Borel compactification of $\mathcal M$.
The curve $C_0$ is mapped to $\Delta/\Gamma$, $C_1$ to $X$ and $C_{\infty}$ to $(D\backslash \Delta)/\Gamma$.
\end{Prop}
\proof It follows from  Theorem 3.5, \cite{m} that $P$ gives an isomorphism between the closure of $\mathcal V$ in the GIT quotient of the space of plane quartics to the Baily Borel compactification of $\mathcal M$. Moreover, strictly semistable points in the closure are mapped to the boundary.
Lemma \ref{0} and Lemma \ref{inf} show that the family $C_{\alpha}$ has a stable reduction in $\alpha=0,\infty$.
In particular it follows from Lemma \ref{0} and \cite{V} that the stable reduction of $C_0$ is the plane quartic with six nodes i.e. the union of four lines. This implies that $X_0$ has period point in $\Delta$ since the extra node gives a $(-2)$ curve in $T$.
By Lemma \ref{inf} $X_{\infty}$ is a $K3$ surface with Picard number $18$, hence its period point is not in $\Delta$ (i.e. the stable reduction of $C_{\infty}$ has only $5$ nodes).
Finally, note that $C_1$ is a strictly semistable quartic, hence it is mapped to $\{X\}$.
\qed

\section{Singular $K3$ surfaces and isogenies}
In this section we study the occurance of singular $K3$ surfaces  in the family and we prove that this is connected to  the existence of isogenies between certain elliptic curves.
\subsection{The curve $B_\alpha$}\label{Balfa}
We consider the genus two curve in Proposition \ref{fi}:
$$B_\alpha:\quad \tau^2=\rho(\rho^4+2\rho^2+\alpha).$$
For further remarks it is convenient to take
$$
\alpha=\beta^{-8}.$$
Now we define $\rho:=\beta^{-2}\rho$, $\tau:=\beta^{-5} \tau$ and the
equation for $B_\alpha$ becomes:
$$
B_\beta:\quad \tau^2=\rho(\rho^4+2\beta^4 \rho^2+1).
$$
It is now easy to see that $B_\beta$ carries the involution
$$
\iota:B_\beta\longrightarrow B_\beta,\qquad (\rho,\tau)\longmapsto
(\rho^{-1},\tau\rho^{-3}).
$$
The quotient by $\iota$ is the elliptic curve:
$$
E_\beta:\quad v^2=u(u^2+4u+2(1+\beta^4))
$$
with quotient map:
$$
f: B_\beta\longrightarrow E_\beta,\qquad (\rho,\tau)\longmapsto
(u,v)=\left(\,\frac{2(1+\beta^4)\rho}{(\rho-1)^2},\;
\frac{2(1+\beta^4)\tau}{(\rho-1)^3}\,\right).
$$
This formula shows that the hyperelliptic involution $(\rho,\tau)\longmapsto (\rho,-\tau)$ on $B_{\beta}$ induces the involution $(u,v)\mapsto (u,-v)$ on $E_\beta$.
\begin{Lem}\label{isog}
The Jacobian of $B_{\beta}$ is isogenous to the
product $E_{\beta}\times E_{\beta}$.
\end{Lem}
\proof The reducibility of the Jacobian of $B_{\beta}$ follows from
Theorem 14.1.1, Ch.14, \cite{Ca2} since it is clear that $B_{\beta}$
is equivalent to a curve of the form:
$$y^2=x(x-1)(x+1)(x-b)(x+b).$$
In fact, the curve $B_{\beta}$ has the automorphism $\iota'$:
$$\iota'(\rho,\tau)=(-\rho,i\tau).$$
This gives another map $f\circ\iota':B_{\beta}\lra E_{\beta}$. Notice
that:
$$H^{1,0}(B_{\beta})=\langle d\rho/\tau, \rho d\rho/\tau\rangle.$$
We have:
$$\begin{array}{ccc}
H^{1,0}(B_{\beta})&=&\langle d\rho/\tau+\rho d\rho/\tau\rangle\oplus\langle d\rho/\tau-\rho d\rho/\tau\rangle,\\
&=&f^*H^{1,0}(E_{\beta})\oplus (f\circ\iota')^*H^{1,0}(E_{\beta}).
\end{array}$$
Hence the Jacobian is isogenous to the product $E_{\beta}\times
E_{\beta}$. \qed
\subsection{The elliptic curve $E$}\label{el}
Notice that we have the isomorphism of genus one curves:
$$
E'=(y^2=x^3-x)\stackrel{\cong}{\longrightarrow} E=(w_1^4=z_1^2-1),$$
$$(x,y)\longmapsto (w_1,z_1)=(y/(\sqrt{2}x),(x+x^{-1})/2)).$$

Moreover, the automorphism of order four on $E$:
$$(z_1,w_1)\longmapsto (z_1,iw_1)$$
is induced by the automorphism on $E'$:
$$(x,y)\longmapsto (x^{-1},iyx^{-2}).$$
We call \emph{standard involution} the automorphism $(x,y)\longmapsto
(x,-y)$ on $E'\cong E$ (sometimes we simply write $p\longmapsto -p$ for
this map).

\subsection{Isogenies}
The construction in the previous section gives a rational map of
degree four from $B_\alpha\times E$ to the quartic surface
$Y_\alpha\subset \Ps 3$.
In coordinates, it is given by:
$$
\Upsilon: B_\alpha\times E\longrightarrow Y_\alpha\subset \Ps 3,$$
$$((\rho,\tau),(z_1,w_1))\longmapsto \left\{
\begin{array}{lcl}
x&=&1,\\
y&=&\rho^2,\\
z&=&\mbox{$\frac{1}{2}$}(\rho^4+2\rho^2+\alpha)z_1+
\mbox{$\frac{1}{2}$}(\rho^4-2\rho^2-\alpha),\\
w&=&\tau w_1/\sqrt{2i}.\\
\end{array}\right.
$$
It can be proved that the image of $B_\beta\times E$ in $\Ps3$ is
the quotient by the order four automorphism
$$
\phi:B_\beta\times E\longrightarrow B_\beta\times E,\qquad
\left((\rho,\tau),(z_1,w_1)\right)\longmapsto
\left((-\rho,i\tau),(z_1,-iw_1)\right).
$$
Note that the square of the automorphism is the product of the
hyperelliptic involution on $B_\beta$ and the standard involution on
$E$. 
\begin{Rem}
The rational map $\Upsilon$ has 9 base points, one of
multiplicity 4 and 8 of multiplicity 2.
\end{Rem}
We now consider the Kummer surface associated to the abelian surface $E_{\beta}\times E$:
$$K_{\beta}=\mbox{Km}(E_{\beta}\times E).$$
From the previous remarks it follows that we have a diagram:
$$\xymatrix{
&(B_{\beta}\times E)/ \phi^2 \ar[ld]_{/\iota} \ar[rd] ^{/\phi} \\
 K_{\beta}&  & X_{\beta} }$$
where the left arrow is the quotient by the involution $\iota$ and the right arrow is the quotient by $\phi$ (composed with birational maps).
We now prove the following
\begin{Thm}\label{kum}
The K3 surface $X_{\beta}$ is singular if and only if $K_{\beta}$ is singular (i.e. $E_{\beta}$ is isogenous
to $E$).
\end{Thm}
\proof Let $\omega=d\rho/\tau\in H^{1,0}(B_{\beta})$,
$\omega_i=dw_1/z_1\in H^{1,0}(E)$. Notice that:
$$(H^1(B_{\beta})\otimes H^1(E))^{\phi}=\langle \omega \otimes \omega_i, \rho\omega\otimes \bar\omega_i, \bar\omega\otimes \bar\omega_i, \overline{\rho\omega}\otimes \omega_i \rangle.$$
Let $\widetilde{B_{\beta}\times E}$ be the blow up of
$B_{\beta}\times E$ along the indeterminacy locus of $\Upsilon$ and
 $\widetilde \Upsilon$ be the map $\widetilde{B_{\beta}\times E}\ra X_{\beta}$ induced by $\Upsilon$.
We have
$$\widetilde{\Upsilon}^*(T_{\beta})\subset (H^1(B_{\beta})\otimes H^1(E))^{\phi}\subset H^2(\widetilde{B_{\beta}\times E})$$
and the first inclusion is an equality for general $\beta\in \Ps
1$. The transcendental lattice $T_{\beta}$ has rank two if the space
$(H^1(B_{\beta},\mathbb Q)\otimes H^1(E,\mathbb Q))^{\phi}$ contains
a cycle of type $(1,1)$. It can be proved by easy computations that
 $H^1(B_{\beta},\mathbb Q)\otimes H^1(E,\mathbb Q)$ is the direct sum
of the eigenspaces (with respect to the eigenvalues $\pm 1$) of the
automorphism $\phi$. Moreover, the involution $\iota$ interchanges
the eigenspaces of $\phi$. This implies that if
$H^1(B_{\beta},\mathbb Q)\otimes H^1(E,\mathbb Q)$ contains a
$(1,1)$ cycle, then the same is true for the positive eigenspace of
$\phi$. Notice that:
$$H^1(B_{\beta})\otimes H^1(E)\cong H^1(B_{\beta})^*\otimes H^1(E)\cong Hom(H^1(B_{\beta}),H^1(E)).$$
Hence we can associate to each element $\omega\in
H^1(B_{\beta})\otimes H^1(E)$ a homomorphism
$\psi_{\omega}:H^1(B_{\beta})\ra H^1(E)$. Moreover $\omega$ is of
type $(1,1)$ iff $\psi_{\omega}$ preserves the Hodge decomposition
i.e. $\psi_{\omega}(H^{1,0}(B_{\beta}))\subset H^{1,0}(E)$ (see
\cite{KS}). By Lemma \ref{isog} $J(B_{\beta})\cong E_{\beta}\times
E_{\beta}$, hence this existence  is equivalent to the existence of
a homomorphism:
$$\psi_{\omega}':H^1(E_{\beta},\mathbb Q)\lra H^1(E,\mathbb Q)$$
preserving the Hodge structure i.e. of an isogeny between
$E_{\beta}$ and $E$.
It is known that the Kummer surface associated to the product of two elliptic curves is singular if and only if the two curves are isogenous with complex multiplication an thus the result follows.
\qed\\

Assume now that $\beta$ is such that there is an isogeny of elliptic
curves:
$$g:E_\beta\lra E.$$
Composing with the quotient map $f$ (see \ref{Balfa}) we have:
$$h:B_\beta\longrightarrow E_\beta\longrightarrow E.$$
Let $\Gamma_h$ be the graph of $h$. By the proof of Theorem \ref{kum} $\Gamma_h$ is the $(1,1)$ cycle in $H^1(B_{\beta})\otimes H^1(E)$ corresponding to $g$.   
\begin{Lem}
The image $\Upsilon(\Gamma_h)$ is a rational curve in $\Ps3$.
\end{Lem}
\proof As observed in \ref{Balfa}, the hyperelliptic involution
$i$ on $B_\beta$ induces the standard involution on $E_\beta$. Since
$g$ is an isogeny (so a homomorphism of groups) the hyperelliptic
involution on $E_{\beta}$ composed with $g$ is the standard
involution on $E$. Thus if $(p,h(p))\in\Gamma_h\cong B_\beta$, then
also $(i(p),h(i(p))=(i(p),-h(p))$ lies in $\Gamma_h$.
This means that the graph $\Gamma_h$ is invariant under $\phi^2$,
therefore the composition
$$
B_\beta\cong\Gamma_h\hookrightarrow B_\beta\times E\lra
Y_{\beta}\subset \Ps3
$$
factors over $B_\beta/i\cong \Ps1$. In particular, the image of the
graph is a rational curve.

\qed

\subsection{A special case}
Theorem 4.5 in \cite{H} predicts that the curve $\Upsilon(\Gamma_h)$ is a ``splitting curve" for $C_{\beta}$ i.e. its inverse image by the cover $\pi_{\beta}$ is the union of four distinct curves.
We prove this in a special example where the isogeny $g:E_\beta\lra E$ is an isomorphism:\\

\noindent\emph{Example:} We consider the curve $B_\beta$ from
section \ref{Balfa} with $\bar\beta^4=7/9$. Then we have:
$$
B_{\bar\beta}:\ \tau^2=\rho(\rho^4+ \mbox{$\frac{14}{9}$}\rho^2+1)
$$
$$
E_{\bar\beta}:\
v^2=u(u+\mbox{$\frac{4}{3}$})(u+\mbox{$\frac{8}{3}$}) .$$ Notice
that, by putting $u=\mbox{$\frac{4}{3}$}x-\mbox{$\frac{4}{3}$}$ we
get an isomorphism with the curve $E':\  y^2=x(x^2-1)$.

We fix the isomorphisms $E_{\bar\beta}\cong E'\cong E$ (the last one
as in section \ref{el}). We denote by $D_{\bar\beta}$ the projection
to $\Ps 2$ of the image of $\Gamma_h$ in $\Ps 3$. Then we have:
\begin{Lem}
The image of $\Gamma_h\subset B_{\bar\beta}\times E$ in $\Ps3$ is a
rational curve of degree six. Moreover, the inverse image of the
curve $D_{\bar\beta}$ splits in four components on the quartic
surface $Y_{\bar\beta}\subset \Ps 3$.
\end{Lem}
\proof

An explicit computation gives that the curve $D_{\bar\beta}\subset
\Ps 2$ is the image of the following map:
$$
\psi:\Ps1_r\lra \Ps2,$$
$$r\longmapsto \left\{
\begin{array}{rcl}
x&=&49(r-1)^2,\\
y&=&63r^2(r-1)^2,\\
z&=&3r^2(48 - 32r + 75r^2 - 54r^3 + 27r^4).
\end{array}\right.
$$
Recall that $Y_{\bar\beta}$ totally ramifies over the plane quartic:
$$
Q:\ (y^2-xz)\cdot y\cdot (\bar\alpha x+2y+z)=0,\ \
\bar\alpha=81/49.
$$
Substituting for $x,y,z$, we get:
$$
\left( -2352( r-1)^2r^2(3 - 2r + 3r^2) \right)\cdot
\left(63r^2(r-1)^2 \right)\cdot \left(3(3 - 2r + 3r^2)^3 \right).
$$
Thus this product is a fourth power in $\C[r]$, hence the 4:1 cover
of the curve splits into 4 components. \qed\\

We remark that the curve $D_{\bar\beta}$ defined in the previous
section defines a 2-section for the elliptic fibration $\mathcal E_{\bar\beta}$
i.e it meets every fiber in two points.
Consider the 2:1 base change:
$$D_{\bar\beta}\lra \Ps1_\lambda$$
given by the projection of the 2-section to the base. The pull-back
$\mathcal E_r$ of the Weierstrass fibration $\mathcal E_{\bar\beta}\lra\Ps1_\lambda$ along this base change has two `new' sections which
are the irreducible components of the pull-back of the 2-section
$D_{\bar\beta}$. The sum of these sections of $\mathcal
E_r\lra \mathbb P^1_r$ actually defines a section of
$\mathcal E_{\bar\beta}\lra \Ps1_\lambda$.
\begin{Lem}
The elliptic fibration $\mathcal E_{\beta}$ with $\beta=\bar\beta$ has a new
section. The inverse image by $\pi_{\beta}$ of the image of this
curve in $\Ps 2$ splits in four components.
\end{Lem}
\proof
The parameter $\lambda$ was defined as $y/x$, (see section \ref{pencil})
hence the base change $D_{\bar\beta}\lra\Ps1$ is defined by:
$$
\lambda=y/x=(9/7)r^2.
$$
On the other hand, the coordinates of the 2-section are polynomials
in $r$, so we need to make a base change with $\sqrt{\lambda}$ or,
equivalently, with $r$. Then the pull-back surface $\mathcal E_r$
has the sections $r\longmapsto x_i(r)$ and $r\longmapsto x_i(-r)$ (here
$x_i=x,y,z,w$ or $u,v$ in the Weierstrass model).

Let $u_1=u(r)$, $u_2=u(-r)$, $v_1=v(r)$, $v_2=v(-r)$. The
coordinates $(u_3,v_3)$ of the sum of the two sections in the
Weierstrass model can be found by using the formula
$$
u_3=(v_2-v_1)^2/(u_2-u_1)^2-u_1-u_2
$$
from \cite{ST} (the $v_3$-coordinate is easy to find from the
Weierstrass equation of $\mathcal E_{\beta}$).
The coordinate $u_3$ is a function of $r^2$ (since $u_1,v_1$ and
$u_2,v_2$ are permuted under $r\mapsto -r$), hence we get a section
of $\mathcal E_{\beta}$. Explicitly, the section of the fibration (note
$\bar\alpha=81/49$):
$$
v^2=u^3-7^{-4}\lambda^3(81 + 98 \lambda + 49 \lambda^2)^2u,
$$
is given by:
$$
u=\frac{(27 + 7\lambda)^2(81 + 98\lambda +
49\lambda^2)}{2^47^{7/2}},$$
$$v=\frac{(81-7\lambda)(27+7\lambda)(81 + 98 \lambda + 49 \lambda^2)^2}
{2^67^{21/4}} .$$
The equation of the corresponding curve in $Y_\alpha$ can be easily
found by using the inverse transformations. The projection of this
curve to $\Ps 2$ is given by (after having replaced $r^2$ by $r$
throughout):
$$
\zeta:\mathbb P^1_r\longrightarrow \mathbb P^2,\ \  r\longmapsto
\left\{
\begin{array}{lcl}
x&=&49(-9 + r)^2,\\
y&=&63r(-9 + r)^2,\\
z&=&9r^2(729 + 94r + 9r^2).
\end{array}\right.
$$
Recall that $Y_\alpha$ totally ramifies over
$$
(y^2-xz)\cdot y\cdot (\alpha x+2y+z),\qquad \bar\alpha=81/49.
$$
Substituting for $x,y,z$, we get:
$$
\left(-9144576 r^3(-9 + r)^2 \right)\cdot \left(63r(-9 + r)^2
\right)\cdot \left( 81 (r+3)^4 \right).
$$
Thus this product is a fourth power in $\C[r]$, hence the 4:1 cover
of the curve splits into four components in $Y_{\alpha}$.
\begin{Rem}
The elliptic fibration $\mathcal E_{\bar\alpha}$ on the $K3$ surface
$X_{\bar\alpha}$ has a fiber of type $III$, one of type $III^*$, two
of type $I_0^*$ and a section. The Shioda-Tate formula (Corollary 1.5, \cite{Shi}) implies that the Picard number of $X_{\bar\alpha}$ is at least 19.
Since the transcendental lattice is a $\Z[i]$ module, the Picard number is an even integer, hence the $K3$ surface is indeed singular.
\end{Rem}

\subsection{The transcendental lattice}

In this
section we prove the following:
\begin{Prop}\label{tra}
The intersection matrix of the transcendental lattice of a singular $K3$ surface $X_{\alpha}$ is of the form:
$$T_n=\left( \begin{array}{cc}
2n&0\\
0&2n\\
\end{array}\right),\ \ n\in\Z,\ n>0.$$
Conversely, if $n$ is a positive integer with $n\not \equiv 2\ (mod\
4)$, then the rank two lattice $T_n$ is the transcendental lattice
of a $K3$ surface $X_{\alpha}$.
\end{Prop}
\proof Let $T_{\alpha}$ be the transcendental lattice of a singular $K3$
surface $X_{\alpha}$ in the family. Notice that $X_{\alpha}$ carries an order four automorphism $\sigma'$ such that the induced isometry $\rho'$ on $H^2(X_{\alpha},\mathbb Z)$ satisfies ${\rho'}^{2}=-id$ on $T_{\alpha}$. Moreover, the transcendental lattice is isomorphic to $\Z[i]$ if we identify $i$ with $\rho'$. It follows (as in the proof of Lemma \ref{iso}) that $T_{\alpha}\cong A_1(-n)^{\oplus 2}$ for some $n$ positive integer.

Conversely, if $n$ is a positive integer, we prove that there
exists $a=(a_1,\dots,a_4)\in\mathbb Z^4$
such that:
$$n=a_1^2+a_2^2-a_3^2-a_4^2,$$
with ${a_1}^2+{a_2}^2>{a_3}^2+{a_4}^2$ and such that the rank two
lattice $\Lambda(a)=\langle a,\rho(a)\rangle$ is primitive in
$\mathbb Z^4$. This is  equivalent to the request that the rank two minors of the matrix:
$$\left(\begin{array}{lccccr}
a_1 & a_2 &a_3 &a_4\\
a_2 & -a_1 & -a_4 &a_3
\end{array}\right)$$
have no common factors. Let $n=2k+1$ be an odd integer, then we can choose: $a_1=(k+1)^2$, $a_3=k^2$ and $a_2=a_4=0$. For $n=2(k+1)$
with odd $k$ we can choose $a_2=1$, $a_4=0$ and $a_1,a_3$ as before.
Assume that $k=2\ell$ is even, then $n=2(2\ell+1)=4\ell+2$. Notice
that:
$$(a_1^2+a_2^2)-(a_3^2+a_4^2)\equiv 2\ (mod\ 4).$$
Since $x^2\in\{0,1\}$ for $x\in \mathbb Z_4$ we have only two
possibilities:
$$1)\ \ (a_1^2,a_2^2,a_3^2,a_4^2)\equiv(1,1,0,0)\ (mod\ 4),$$
$$2)\ \ (a_1^2,a_2^2,a_3^2,a_4^2)\equiv(0,0,1,1)\ (mod\ 4).$$
In case $1)$ we have that $a_1,a_2$ are odd and $a_3,a_4$ is even.
Hence we immediately get that $\Lambda(a)$ is not primitive (all
minors are even integers). The second case is analogous.

We now assume that $n\not\equiv 2\ (mod\ 4)$ and  we choose
$a_1,\dots,a_4$ as before. Define $z(a)\in T\otimes\C$ by (with
respect to the usual basis):
$$z(a)=(a_1+ia_2, a_2-ia_1,a_3+ia_4, a_4-ia_3).$$
We consider the sublattices of $L$:
$$N(a)={z(a)}^{\perp}\cap L_{K3},\ \ T(a)=N(a)^{\perp}.$$
Notice that $T(a)=<a,\rho(a)>$ with intersection matrix  given by:
$$T(a)\cong\left(\begin{array}{cc}
2({a_1}^2+{a_2}^2-{a_3}^2-{a_4}^2)& 0\\
0&2({a_1}^2+{a_2}^2-{a_3}^2-{a_4}^2)\\
\end{array}\right).$$

By the surjectivity of the period map,
for every $a\in \mathbb Z^4$ as above there exists a marked $K3$
surface $X(a)$ with period point $z(a)$. In particular, the
transcendental lattice of $X(a)$ is $T(a)$. \qed\\

We now give some examples of $K3$ surfaces with transcendental lattice isomorphic to $T_n$:
\begin{itemize}
\item[a)] $n=1$ for \emph{Vinberg's K3 surface} (see \cite{V}),
\item[b)] $n=3$ for the $K3$ surface described in \cite{KOZ},
\item[c)] $n=4$ for the \emph{Fermat quartic} (see \cite{O1}),
\item[d)] $n=7$ for the \emph{Klein quartic} (see \cite{OZ}).
\end{itemize}
\begin{Lem}\label{m}
Let $E\cong \C/\Z+i\Z$ and $E'\cong \C/\Z+mi\Z$, then the Kummer surface $\mbox{Km}(E\times E')$  has transcendental lattice of the form $T_{2m}$, $m\in \Z$, $m>0$. 
\end{Lem}
\proof It follows easily from the proof of Theorem 4, \cite{SI}.\qed
 \begin{Lem}
The family $\{X_{\alpha}\}$ contains the $K3$ surfaces  a), b), c), d) and all Kummer surfaces in Lemma \ref{m} with even $m$. The surface $\mbox{Km}(E\times E)$ is not in the family.
\end{Lem}
\proof The first assertion is a corollary of Proposition \ref{tra}. The transcendental lattice of $X=\mbox{Km}(E\times E)$ is isomorphic to $T_2$ (see also \cite{KK}), so Proposition \ref{sing} can not be applied. Assume that $X=X_{\alpha}$, $\alpha\in \Ps 1$. Notice that $X$ can not correspond to the fibers $\alpha=0,\infty$, since it is singular and it is not isomorphic to Vinberg's $K3$ surface. Hence the elliptic fibration $\mathcal E_{\alpha}$ on $X$ has the same configuration of singular fibers of the general case i.e. of Lemma \ref{sing}. Since the rank of the Picard lattice is $20$, it follows that the Mordell-Weil
group modulo torsion is isomorphic to $\Z^{\oplus 2}$. However, the table in \cite{Nis}, shows that there exists no elliptic fibration on $X$ with these properties. \qed

\bibliographystyle{plain}

\end{document}